\documentclass[a4paper, 11pt]{amsart}

\usepackage[mathscr]{eucal}
\usepackage{amsmath}
\usepackage{amssymb}
\usepackage{amsfonts}
\usepackage{amsthm}

\usepackage{cases}

\usepackage[mathscr]{eucal}
\usepackage{amsmath}
\usepackage{amssymb}
\usepackage{amsfonts}
\usepackage{amsthm}
\usepackage[dvips]{graphics, color}

\theoremstyle{plain}
\newtheorem{theorem}{Theorem}[section]
\newtheorem{proposition}{Proposition}[section]

\theoremstyle{remark}

\newtheorem{remark}{Remark}[section]

\numberwithin{equation}{section}

\def\<{\left<} \def\>{\right>}

\def\proof{\noindent{\it Proof. }}
\def\be{\begin{equation} }
\def\ee{\end{equation} }
\def\qed{\ifhmode\unskip\nobreak\fi\ifmmode\ifinner\else\hskip5pt
\fi\fi\hbox{\hskip5 pt \vrule width4 pt height6 pt depth1.5 pt \hskip1pt }}

\begin{document}

\title[]{Surfaces in Euclidean 3-space with Maslovian normal bundles}
\author[]{Toru Sasahara}
\address{Division of Mathematics, Center for Liberal Arts and Sciences, 
Hachinohe Institute of Technology, 
Hachinohe, Aomori, 031-8501, Japan}
\email{sasahara@hi-tech.ac.jp}

\date{}

\begin{abstract}
{\footnotesize
We prove that a surface in Euclidean $3$-space has Maslovian normal bundle if and only if 
it is  a part of a round sphere,  a circular cylinder, or a circular cone.
}
\end{abstract}


 \maketitle

 \section{Introduction}
 An $n$-dimensional submanifold $M$ in a K\"{a}hler $n$-manifold $N$
is called
Lagrangian if the complex structure $J$ of $N$  interchanges the tangent and normal spaces of $M$.
The study of Lagrangian submanifolds from the viewpoint of Riemannian geometry has
been an  active field over the last half century.


The normal bundle of a submanifold in Euclidean $n$-space ${\mathbb R}^n$ 
 can be naturally  immersed 
in ${\mathbb C}^n$ as a Lagrangian submanifold (see \cite{hl}). 
Thus, it is natural to study submanifolds in  ${\mathbb R}^n$  by
imposing some special conditions on the normal bundles.
Harvey and Lawson \cite{hl} showed that a submanifold in ${\mathbb R}^n$ 
has minimal normal bundle if and only if for each normal vector, the set of  eigenvalues of its shape operator 
is invariant under the multiplication by $-1$.
This implies that a surface  in ${\mathbb R}^n$ 
has minimal normal bundle if and only if it is minimal.
Sakaki \cite{sakaki} proved that  a surface in ${\mathbb R}^3$ has Hamiltonian stationary normal bundle if and only if it is either minimal, a part of  a round sphere, or a part of 
a cone with vertex angle $\pi/2$.
 It  was  proved in \cite{sasa} that a surface in ${\mathbb R}^3$ has tangentially biharmonic normal bundle if and only if it is either minimal, a part of a round sphere, or a part of a circular cylinder.

For  a Lagrangian submanifold, 
 the dual form of $JH$, where $H$ is the mean curvature vector field, 
 is the Maslov form (up to a constant).
A Lagrangian submanifold is called {\it Maslovian} if 
$H$ vanishes  nowhere  and 
$JH$ is a principal direction
of $A_H$, where $A_{H}$ is the shape operator with respect to $H$ (see \cite{chen}).
Many important Lagrangian submanifolds satisfy the Maslovian condition (see, for example,
\cite{cas,  chd, chv, ros}).



 In this paper, we determine  all surfaces in ${\mathbb R}^3$ with Maslovian normal bundles as follows:
 \begin{theorem}\label{main}
 A surface in ${\mathbb R}^3$ has Maslovian normal bundle if and only if 
it is a part of a round sphere, a circular cylinder, or a circular cone.
 \end{theorem}
 
\section{Preliminaries}
Let $M$ be  
a  submanifold of 
a Riemannian manifold $\tilde M$ and $\iota$ its immersion.
We identify a point $x\in M$ with  $\iota(x)$ and a tangent vector $X\in T_xM$ 
with $\iota_{*}(X)$.
 We denote by $\nabla$ and
$\tilde\nabla$ the Levi-Civita connections on $M$ and 
$\tilde M$, respectively. The
formulas of Gauss and Weingarten are given respectively by
\begin{equation}\label{gw}
\tilde \nabla_XY= \nabla_XY+h(X,Y), \quad  \tilde\nabla_X \xi = -A_{\xi}X+D_X\xi,
\end{equation}
 for tangent vector fields $X$, $Y$ and a normal vector field $\xi$, 
where $h,A$ and $D$ are the second fundamental
form, the shape operator and the normal
connection.
In this paper, the mean curvature vector field $H$ is defined as 
$H={\rm trace}\hskip3pt h$.


If $M$ is a 
 hypersurface of ${\mathbb R}^n$, then
the Gauss and Codazzi equations are given respectively by
\begin{gather}
R(X, Y)Z=\<AY, Z\>AX-\<AX, Z\>AY, \label{gau}\\
(\nabla_XA)Y=(\nabla_YA)X, \label{cod}
\end{gather}
where $R$ is the curvature tensor of $M$ and  $A$ is the shape operator with respect to 
the unit normal vector field. 
\section{Normal bundles of surfaces in ${\mathbb R}^3$} 
Let $M$ be a surface in ${\mathbb R}^3$.
The normal bundle $T^{\perp}M$ of $M$ is naturally immersed in ${\mathbb R}^3
\times{\mathbb R}^3$ by the immersion 
$f(\xi_x):=(x, \xi_x)$, which is  expressed as
\be\label{nb}
f(x, t)=(x, tN)
\ee
for $t\in{\mathbb R}$ and the unit normal vector field $N$ along $x$.
We equip $T^{\perp}M$ with the metric induced by $f$. 


We choose   a local orthonormal frame $\{e_1, e_2\}$ on an open subset  $U$ of $M$ such that
\begin{equation}\label{pri}
Ae_1=ae_1,\quad
Ae_2=be_2
\end{equation}
for some functions $a$ and $b$.
Put $\<\nabla_{e_i}e_j, e_k\>=\omega_j^k(e_i)$ for $i, j ,k\in\{1, 2\}$. 
Note that $\omega_1^2=-\omega_2^1$.
The Codazzi equation (\ref{cod}) yields
\begin{equation}\label{coda}
e_1b=(a-b)\omega_1^2(e_2), \quad e_2a=(b-a)\omega_2^1(e_1).
\end{equation}
We define the following tangent vector fields on $U\times{\mathbb R}\subset T^{\perp}M$:
\begin{equation}
\begin{split}\label{e123}
&\tilde e_1=(1+t^2a^2)^{-\frac{1}{2}}e_1, \\
&\tilde e_2=(1+t^2b^2)^{-\frac{1}{2}}e_2, \\
&\tilde e_3=\frac{\partial}{\partial t}.
\end{split}
\end{equation}
From (\ref{gw}), (\ref{nb}) and (\ref{pri}), it follows   that
\begin{equation}
\begin{split}\label{fe123}
&f_{*}(\tilde e_1)=(1+t^2a^2)^{-\frac{1}{2}}(e_1, -tae_1),\\
&f_{*}(\tilde e_2)=(1+t^2b^2)^{-\frac{1}{2}}(e_2, -tbe_2),\\
&f_{*}(\tilde e_3)=(0, N).
\end{split}
\end{equation}
Thus, $\{\tilde e_1, \tilde e_2, \tilde e_3\}$ is an orthonormal frame on
 $U\times{\mathbb R}$.
 
Let $J$ be the complex structure on ${\mathbb C}^3={\mathbb R}^3\times{\mathbb R}^3$
by $J(X, Y):=(-Y, X)$.
We define the following vector fields along $f$:
\begin{equation}
\begin{split}\label{e456}
&e_4:=Jf_{*}(\tilde e_1)=(1+t^2a^2)^{-\frac{1}{2}}(tae_1, e_1),\\
&e_5:=Jf_{*}(\tilde e_2)=(1+t^2b^2)^{-\frac{1}{2}}(tbe_2, e_2),\\
&e_6:=Jf_{*}(\tilde e_3)=(-N, 0).
\end{split}
\end{equation}
Then $\{e_4, e_5, e_6\}$ is a normal orthonormal frame. This implies that
$T^{\perp}M$ is a Lagrangian submanifold of ${\mathbb C}^3$.

Put
$h^{\alpha}_{ij}=\langle{\tilde e_i}(f_{*}(\tilde e_j)), e_{\alpha}\rangle$
for  $1\leq i, j\leq 3$,  $4\leq\alpha\leq 6$. 
It follows from (\ref{gw})   that
 the mean curvature vector field $H$ of $T^{\perp}M$ in
$\mathbb{C}^3$ is given by 
$H=\sum_{\alpha=4}^{6}\sum_{i=1}^3h_{ii}^{\alpha}e_\alpha$.
From (\ref{coda})-(\ref{e456}), we have
\begin{equation}
\begin{split}\label{second}
h^4_{11}&=-t(1+t^2a^2)^{-\frac{3}{2}}e_1a,\\
h^4_{22}&=-t(1+t^2a^2)^{-\frac{1}{2}}(1+t^2b^2)^{-1}(b-a)\omega_2^1(e_2),\\
&=-t(1+t^2a^2)^{-\frac{1}{2}}(1+t^2b^2)^{-1}e_1b,\\
h^5_{11}&=-t(1+t^2a^2)^{-1}(1+t^2b^2)^{-\frac{1}{2}}(a-b)\omega_1^2(e_1),\\
&=-t(1+t^2a^2)^{-1}(1+t^2b^2)^{-\frac{1}{2}}e_2a,\\
h^5_{22}&=-t(1+t^2b^2)^{-\frac{3}{2}}e_2b,\\
h^{6}_{11}&=-a(1+t^2a^2)^{-1},\\
h^{6}_{22}&=-b(1+t^2b^2)^{-1},\\
h^{4}_{33}&=h^5_{33}=h^6_{33}=0.
\end{split}
\end{equation}
Using (\ref{e456}) and (\ref{second}), we obtain (see \cite{sakaki} and \cite{sasa})
\be
\begin{split}\label{HJH}
H&=-(Pt^2ae_1+Qt^2be_2-RN, Pte_1+Qte_2),\\
JH&=tPe_1+tQe_2+R\frac{\partial}{\partial t},
\end{split}
\ee
where $P$, $Q$ and $R$ are given by 
\begin{equation}
\begin{split}\label{PQR}
&P=(1+t^2a^2)^{-2}e_1a+(1+t^2a^2)^{-1}(1+t^2b^2)^{-1}e_1b,\\
&Q=(1+t^2a^2)^{-1}(1+t^2b^2)^{-1}e_2a+(1+t^2b^2)^{-2}e_2b,\\
&R=a(1+t^2a^2)^{-1}+b(1+t^2b^2)^{-1}.
\end{split}
\end{equation}

From (\ref{HJH}) and (\ref{PQR}),  we obtain the following (cf. [4, III. Th.3.11, Pop.2.17]):
\begin{proposition}\label{pro}
 A surface $M$ in ${\mathbb R}^3$ is minimal if and only if 
$T^{\perp}M$ is a minimal submanifold of ${\mathbb C}^3$.
\end{proposition}
The following two theorems are generalizations of Proposition \ref{pro}.
\begin{theorem}[\cite{sakaki}]\label{sakaki}
Let M be a surface in ${\mathbb R}^3$. Then $T^{\perp}M$ is Hamiltonian stationary
if and only if  $M$ is either minimal, a part of  a round sphere, or a part of 
a cone with vertex angle $\pi/2$.
\end{theorem}
\begin{theorem}[\cite{sasa}]
A surface in ${\mathbb R}^3$  has tangentially biharmonic normal bundle
 if and only if it is either minimal, a part of 
a round sphere, or
a part of a circular cylinder.
\end{theorem}

\begin{remark}
 The notion of tangentially biharmonic submanifolds was introduced by the author in \cite{sasa}.
 This notion agrees with  that  of biconservative submanifolds introduced in \cite{ca}.
Many interesting results  on this subject have been obtained  in the last decade
 (see, for example, \cite{fu, nis, sasa2, tur}
 and references therein).
\end{remark}

\section{Proof of Theorem \ref{main}}
\proof
Let $M$ be a surface in 
${\mathbb R}^3$. 
We denote by $A$ the shape operator of $T^{\perp}M$ in ${\mathbb C}^3$.
Note that $A_{H}JH$ is the tangential part of $-JH(H)$.
For $1\leq i<j\leq 3$, we put 
\be
F_{ij}=\<JH(H), f_{*}(\tilde e_i)\>\<JH, \tilde e_j\>
-\<JH(H), f_{*}(\tilde e_j)\>\<JH, \tilde e_i\>.
\ee
Then, $T^{\perp}M$ is a
Maslovian
Lagrangian submanifold 
in ${\mathbb C}^3$, that is, $A_H(JH)$ is parallel to $JH$, 
if and only if  
\be 
F_{12}=F_{13}=F_{23}=0.\label{mas}
\ee

We shall compute $\<JH(H), f_{*}(\tilde e_i)\>$ and $\<JH, \tilde e_i\>$
for $i=1, 2, 3$.
By (\ref{HJH}), we have
\be
\begin{split}
JH(H)=\biggl(& tP\Bigl[-(e_1P)t^2ae_1-Pt^2(e_1a)e_1-Pt^2a\{\omega_1^2(e_1)e_2+aN\}\\
& -(e_1Q)t^2be_2-Qt^2(e_1b)e_2-Qt^2b\omega_2^1(e_1)e_1+(e_1R)N-aRe_1\Bigr]\\
& +tQ\Bigl[-(e_2P)t^2ae_1-Pt^2(e_2a)e_1-Pt^2a\omega_1^2(e_2)e_2-(e_2Q)t^2be_2\\
&-Qt^2(e_2b)e_2-Qt^2b\{\omega_2^1(e_2)e_1+bN\}+(e_2R)N-bRe_2\Bigr]\\
&+R\Bigl[-\frac{\partial P}{\partial t}t^2ae_1-2tPae_1-\frac{\partial Q}{\partial t}t^2be_2
-2tQbe_2+\frac{\partial R}{\partial t}N\Bigr],\\
& tP\Bigl[-(e_1P)te_1-Pt\{\omega_1^2(e_1)e_2+aN\}-(e_1Q)te_2-Qt\omega_2^1(e_1)e_1\Bigr]\\
&+tQ\Bigl[-(e_2P)te_1-Pt\omega_1^2(e_2)e_2-(e_2Q)te_2-Qt\{\omega_2^1(e_2)e_1+bN\}\Bigr]\\
&+R\Bigl[-\frac{\partial P}{\partial t}te_1-Pe_1-\frac{\partial Q}{\partial t}te_2-Qe_2\Bigr]\biggr).\label{H1}
\end{split}
\ee
Using (\ref{H1}), we have
\be 
\begin{split}
\<JH(H), f_{*}(\tilde e_1)\>=&(1+t^2a^2)^{-\frac{1}{2}}\Bigl[-t^3P(e_1P)a-P^2t^3e_1a-PQt^3b\omega_2^1(e_1)\\
&-atPR-t^3Q(e_2P)a-t^3PQ(e_2a)-t^3Q^2b\omega_2^1(e_2)\\
&-\frac{\partial P}{\partial t}Rt^2a-2tRPa+t^3P(e_1P)a
+t^3PQa\omega_2^1(e_1)\\
&+t^3Q(e_2P)a+t^3Q^2\omega_2^1(e_2)a+R\frac{\partial P}{\partial t}t^2a+tPRa\Bigl]\\
=&(1+t^2a^2)^{-\frac{1}{2}}\Bigl[-P^2t^3e_1a-PQt^3b\omega_2^1(e_1)-t^3PQe_2a\\
&-t^3Q^2b\omega_2^1(e_2)
-2tPRa+t^3PQa\omega_2^1(e_1)+t^3Q^2\omega_2^1(e_2)a\Bigr],\label{H2}
\end{split}
\ee
\be 
\begin{split}
\<JH(H), f_{*}(\tilde e_2)\>=&(1+t^2b^2)^{-\frac{1}{2}}
\Bigl[-t^3P^2a\omega_1^2(e_1)-t^3P(e_1Q)b-t^3PQ(e_1b)\\
&-t^3PQa\omega_1^2(e_2)-t^3Q(e_2Q)b-t^3Q^2(e_2b)-tQRb\\
&-R\frac{\partial Q}{\partial t}t^2b-2tRQb
+t^3P^2b\omega_1^2(e_1)
+t^3P(e_1Q)b\\
&+t^3PQ\omega_1^2(e_2)b+t^3Q(e_2Q)b+t^2\frac{\partial Q}{\partial t}Rb+tRQb\Bigr]\\
=&(1+t^2b^2)^{-\frac{1}{2}}
\Bigl[-t^3P^2a\omega_1^2(e_1)-t^3PQ(e_1b)-t^3PQa\omega_1^2(e_2)\\
&-t^3Q^2(e_2b)-2tRQb+t^3P^2b\omega_1^2(e_1)+t^3PQ\omega_1^2(e_2)b\Bigr].\label{H3}
\end{split}
\ee
Applying (\ref{coda}) to (\ref{H2}) and (\ref{H3}) gives
\be 
\begin{split}
\<JH(H), f_{*}(\tilde e_1)\>&=(1+t^2a^2)^{-\frac{1}{2}}(-t^3P^2e_1a-2t^3PQe_2a-t^3Q^2e_1b-2tPRa),\\
\<JH(H), f_{*}(\tilde e_2)\>&=(1+t^2b^2)^{-\frac{1}{2}}(-t^3P^2e_2a-2t^3PQe_1b
-t^3Q^2e_2b-2tRQb). \label{H4}
\end{split}
\ee
From (\ref{H1}), we easily obtain  
\be 
\<JH(H), f_{*}(\tilde e_3)\>=-t^2P^2a-t^2Q^2b. \label{H5}
\ee
It follows from (\ref{e123}), (\ref{HJH}) and (\ref{PQR}) that
\be
\begin{split}
&\<JH,  \tilde e_1\>=t\Bigl[(1+t^2a^2)^{-\frac{3}{2}}e_1a
+(1+t^2a^2)^{-\frac{1}{2}}(1+t^2b^2)^{-1}e_1b\Bigr],\\
&\<JH, \tilde e_2\>=t\Bigl[(1+t^2a^2)^{-1}(1+t^2b^2)^{-\frac{1}{2}}e_2a
+(1+t^2b^2)^{-\frac{3}{2}}e_2b\Bigr],\\
&\<JH,  \tilde e_3\>=a(1+t^2a^2)^{-1}+b(1+t^2b^2)^{-1}.\label{JH1}
\end{split}
\ee

{\bf Case (I).} $M$ is isoparametric. In this case,
 $M$ is either a part of a round sphere or 
a part of a circular cylinder. By (\ref{H4}), (\ref{H5}) and (\ref{JH1}),
we see that  (\ref{mas}) is clearly satisfied. Hence, $T^{\perp}M$ is Maslovian.

{\bf Case (II).} $M$ is non-isoparametric. 
Using (\ref{H4}), (\ref{H5}) and (\ref{JH1}), we have
\be 
\begin{split}
(1+t^2a^2)^{\frac{1}{2}}F_{13}=&a(1+t^2a^2)^{-1}(-P^2t^3e_1a-2t^3PQe_2a-t^3Q^2e_1b-2tPRa)\\
&+b(1+t^2b^2)^{-1}(-P^2t^3e_1a-2t^3PQe_2a-t^3Q^2e_1b-2tPRa)\\
&+t^3(P^2a+Q^2b)\bigl[(1+t^2a^2)^{-1}e_1a+(1+t^2b^2)^{-1}e_1b\bigr]\\
=&a(1+t^2a^2)^{-1}(-2t^3PQe_2a-t^3Q^2e_1b-2tPRa)\\
&+b(1+t^2b^2)^{-1}(-P^2t^3e_1a-2t^3PQe_2a-2tPRa)\\
&+t^3\bigl[(1+t^2b^2)^{-1}P^2ae_1b+(1+t^2a^2)^{-1}Q^2be_1a\bigr]\\
=&\phi_1(x, t)t^3+\phi_2(x,t)t,
\end{split}\label{F13}
\ee
where $\phi_1(x, t)$ and $\phi_2(x, t)$ are functions on $U\times {\mathbb R}$ given by
\begin{align*}
\phi_1(x, t)=&-2PQRe_2a+(1+t^2a^2)^{-1}Q^2(be_1a-ae_1b)
+(1+t^2b^2)^{-1}P^2(ae_1b-be_1a),\\
\phi_2(x, t)=&-2aPR^2.
\end{align*}
In the same way as above, we have
\be 
\begin{split}
(1+t^2b^2)^{\frac{1}{2}}F_{23}=&a(1+t^2a^2)^{-1}(-t^3P^2e_2a-2t^3PQe_1b-t^3Q^2e_2b-2tRQb)\\
&+b(1+t^2b^2)^{-1}(-t^3P^2e_2a-2t^3PQe_1b-t^3Q^2e_2b-2tRQb)\\
&+t^3(P^2a+Q^2b)\bigl[(1+t^2a^2)^{-1}e_2a+(1+t^2b^2)^{-1}e_2b\bigr],\\
=&a(1+t^2a^2)^{-1}(-2t^3PQe_1b-t^3Q^2e_2b-2tRQb)\\
&+b(1+t^2b^2)^{-1}(-t^3P^2e_2a-2t^3PQe_1b-2tRQb)\\
&+t^3\bigl[(1+t^2b^2)^{-1}P^2ae_2b+(1+t^2a^2)^{-1}Q^2be_2a\bigr]\\
=&\psi_1(x, t)t^3+\psi_2(x, t)t,
\end{split}\label{F23}
\ee
where  $\psi_1(x, t)$ and $\psi_2(x, t)$ are functions on $U\times {\mathbb R}$ given by
\begin{align*}
\psi_1(x, t)=&-2PQRe_1b+(1+t^2a^2)^{-1}Q^2(be_2a-ae_2b)+(1+t^2b^2)^{-1}P^2(ae_2b-be_2a),\\
\psi_2(x, t)=&-2bQR^2.
\end{align*}
We substitute (\ref{PQR}) into the right-hand sides of (\ref{F13}) and (\ref{F23}).
Then, multiplying $(1+t^2a^2)^4(1+t^2b^2)^4$ on both sides of
(\ref{F13}) and (\ref{F23}), we find
\begin{align} 
&(1+t^2a^2)^{\frac{9}{2}}(1+t^2b^2)^4F_{13}=\sum_{i=1}^{5}f_{2i-1}(
a, b, e_1a, e_1b, e_2a, e_2b)t^{2i-1}, \label{FF13}\\
&(1+t^2a^2)^4(1+t^2b^2)^{\frac{9}{2}}F_{23}=\sum_{i=1}^{5}g_{2i-1}(a, b, e_1a, e_1b, e_2a, e_2b)t^{2i-1} \label{FF23}
\end{align}
for some polynomials $f_{2i-1}$ and $g_{2i-1}$ in $a$, $b$, $e_1a$, $e_1b$, $e_2a$
and  $e_2b$.
It is not difficult to see that $f_1$ and $g_1$ coincide with $\phi_2(x, 0)$ and $\psi_2(x, 0)$, respectively.
Thus, we obtain 
\be
\begin{split}
f_1&=-2a(a+b)^2(e_1a+e_1b),\\
g_1&=-2b(a+b)^2(e_2a+e_2b).
\end{split}\label{fg1}
\ee

{\bf Case (II.1).} $ab\ne 0$. If $T^{\perp}M$ is Maslovian, then  (\ref{mas}) implies that
 (\ref{FF13}) and (\ref{FF23}) are identically zero, and hence $f_{2i-1}=g_{2i-1}=0$ for $i=1, 
 2, 3, 4, 5$. 
Thus, (\ref{fg1}) yields
\be 
e_1a+e_1b=e_2a+e_2b=0. \label{al} 
\ee
We substitute $e_1b=-e_1a$ and $e_2b=-e_2b$
 into (\ref{F13}) and (\ref{F23}).
 The constant terms in  $(1+t^2a^2)^4(1+t^2b^2)^4\phi_1(x, t)$ and 
 $(1+t^2a^2)^4(1+t^2b^2)^4\psi_1(x, t)$ coincide with $\phi_1(x, 0)$ and $\psi_1(x, 0)$,
 respectively.  If $t=0$, then $P=Q=0$, and hence $\phi_1(x, 0)=\psi_1(x, 0)=0$.
 This implies that $f_3$ and $g_3$ coincide with the coefficients of $t^2$ in $(1+t^2a^2)^4(1+t^2b^2)^4\phi_2(x, t)$ and $(1+t^2a^2)^4(1+t^2b^2)^4\psi_2(x, t)$,
 respectively. Thus, by a straightforward computation, we find that  $f_3$ and $g_3$ can be
reduced to the following simple forms:
\be 
\begin{split}
&f_3=2a(a-b)(a+b)^3e_1a, \\
&g_3=2b(a-b)(a+b)^3e_2a.
\end{split}\notag
\ee
Therefore, we have  $e_1a=e_2a=0$ because of $a\ne b$. Combining this with (\ref{al})
yields
that   $a$ and $b$ are constant, which is a contradiction.
Consequently, in this case $T^{\perp}M$ can not be Maslovian.

{\bf Case (II.2).} $ab=0$.
 We assume that 
$a=0$ and $b$ is not constant. 
The relation (\ref{F23}) becomes
\be 
(1+t^2b^2)^{\frac{1}{2}}F_{23}=-2b(1+t^2b^2)^{-4}\bigl[t^3(e_1b)^2e_2b+t(e_2b)b^2\bigr]. \label{f23}
\ee
If $T^{\perp}M$ is Maslovian, then (\ref{f23}) is identically zero, and hence $e_2b=0$, which leads to $Q=0$. 
Hence, from (\ref{H4}) and (\ref{F13}), we  find  that
$F_{12}=0$ and $F_{13}=0$ are automatically satisfied.

Using (\ref{coda}), we have $[e_1, e_2/b]=0$. Thus,  
there exist  local coordinates
$\{t_1, t_2\}$ such that
\be
e_1=\frac{\partial}{\partial t_1}, \quad e_2=b\frac{\partial}{\partial t_2}.\nonumber
\ee
Since $e_2b=0$, we have $b=b(t_1)$.
By the Gauss equation (\ref{gau}), we see  that $M$ is flat. Therefore,  
we get
\be 
b(t_1)=\frac{1}{rt_1+c} \nonumber
\ee for some constants $r$ and  $c$.
Since $b$ is not constant, we have $r\ne 0$. 
After  the coordinate transformation: 
\be
u=t_1+\frac{c}{r}, \quad v=rt_2,
\ee  the metric tensor $g$ and the
second fundamental form $h$ take the following forms:
\be 
g=du^2+u^2dv^2, \quad h=\frac{u}{r}dv^2. \label{gh}
\ee

The circular cone given by
\be 
x(u, v)=\frac{u}{\sqrt{r^2+1}}\Bigl(r\cos\Bigl(\frac{\sqrt{r^2+1}}{r}v\Bigr),
r\sin\Bigl(\frac{\sqrt{r^2+1}}{r}v\Bigr), 1\biggr) \label{cone}
\ee
has the metric tensor and the second fundamental form described in 
(\ref{gh}).
By the fundamental theorem in the theory of surfaces, $M$ is congruent to a part of  (\ref{cone}).

Conversely, if $M$ is parametrized by (\ref{cone}), then  substituting  
$e_1=\frac{\partial}{\partial u}$, $e_2=\frac{1}{u}\frac{\partial}{\partial v}$, $a=0$
 and $b=1/(ru)$ into (\ref{H4}) and (\ref{H5}),  we find 
  that  (\ref{mas}) is satisfied, and hence $T^{\perp}M$ is Maslovian. 
The proof is finished. \qed

\begin{remark}
From (\ref{H4}) and (\ref{H5}), we see that 
every normal bundle in Theorem \ref{main} satisfies $A_{JH}H=0$.
The normal bundle of  the cone  (\ref{cone}) with $r=1$
 is Hamiltonian
 stationary (see Theorem \ref{sakaki}). 
\end{remark}


\begin{thebibliography}{99}
\bibitem{ca}
R. Caddeo, S. Montaldo, C. Oniciuc and P. Piu,
{Surfaces in three-dimensional space forms with
divergence-free stress-bienergy tensor,}
Ann. Mat. Pura Appl. {\bf 193} (2014), 529-550.



\bibitem{cas}
I. Castro and F. Urbano, {\it Twistor holomorphic Lagrangian surfaces in complex
projective and hyperbolic planes}, Ann. Global Anal. Geom. {\bf 13} (1995), 59-67. 








\bibitem{chen}
B. Y. Chen and O. J. Garay,
{\it Maslovian Lagrangian isometric immersions of real space forms
into complex space forms,}
Japan. J. Math. {\bf 30} (2004), 227-281.

\bibitem{chd}
B. Y. Chen, F. Dillen, J. Van der Veken and L. Vrancken,
{\it Curvature inequalities for Lagrangian submanifolds: the final solution,}
Differential Geom. Appl. {\bf 31} (2013), 808-819.


\bibitem{chv}
B. Y. Chen, L. Vrancken and X. Wang,
{\it Lagrangian submanifolds in complex space forms satisfying equality in the optimal
inequality involving $\delta(2, \ldots, 2)$,}
Beitr. Algebra Geom. {\bf 62} (2021), 251-264.

\bibitem{fu}
Y. Fu and N. C. Turgay, {\it Complete classification of biconservative hypersurfaces with 
diagonalizable shape operator in the Minkowski $4$-space,} 
Internat. J. Math. {\bf 27} (2016), no.5,  1650041, 17 pp.


\bibitem{hl}
R. Harvey and H. B. Lawson, {\it Calibrated geometries,}
Acta Math. {\bf 148} (1982), 47-157.

\bibitem{nis}
S. Nistor and C. Oniciuc, {\it Complete biconservative surfaces in the hyperbolic space
${\mathbb H}^3$,} Nonlinear Anal. {\bf 198} (2020), 111860, 29 pp.


\bibitem{ros}
A. Ros and F. Urbano, {\it Lagrangian submanifolds of  ${\mathbb C}^n$
with conformal Maslov form and the Whitney sphere,}
J. Math. Soc. Japan {\bf 50} (1998), 203-226.

\bibitem{sakaki}
M. Sakaki, {\it Hamiltonian stationary normal bundles of surfaces in ${\bf R}^3$},
Proc. Amer. Math. Soc. {\bf 127}(1999), 1509-1515.


\bibitem{sasa}
T. Sasahara, {\it Surfaces in Euclidean $3$-space whose normal bundles are
tangentially biharmonic}, Arch. Math. (Basel) {\bf 99} (2012), 281-287.

\bibitem{sasa2}
T. Sasahara, {\it Tangentially biharmonic Lagrangian $H$-umbilical
submanifolds in complex space forms,}
Abh. Math. Semin. Univ. Hamb. {\bf 85} (2015), 107-123.

\bibitem{tur}
N. C. Turgay and A.  Upadhyay, {\it On biconservative hypersurfaces in $4$-dimensional
Riemannian space forms,} Math. Nachr. {\bf 292} (2019), 905-921.





\end{thebibliography}
\end{document}